\documentclass[12pt,thmsa]{article}
\usepackage{sw20lart}

\input{tcilatex}
\begin{document}

\author{Alan Horwitz}
\title{Even Compositions of Entire Functions and Related Matters }

\section{Introduction}

Let $f(z)$\ and $g(z)$\ be entire functions with their composition $f\circ g$%
\ even. What can be said about $f$\ and $g$\ ? Also, given certain
conditions on $f$\ and/or $g,$\ can $f\circ g$\ be even ? There are numerous
examples which show that $f\circ g$\ can be even, while neither $f$\ nor $%
g\; $is even. In \cite{HR} we proved that if $f$\ and $g$\ are polynomials
with $g(0)=0$\ and $g$\ neither even nor odd, then $f\circ g$\ cannot be
even. It was also shown that if $f$\ and $g$\ are polynomials, $f$ not
constant, $g(0)=0,$\ and $f\circ g$\ is even, then $f$\ or $g$\ must be
even. These results do not hold for entire functions in general, or for
rational functions in general.

Call an entire function $g(z)$\ right pseudo even(RPE) if there is a
non-constant entire function $f(z)$\ such that $f\circ g$\ is even. Of
course, every even or odd entire function is RPE, but there are entire
functions $g$\ which are RPE, while $g$\ (or any derivative of $g)$\ is
neither even nor odd. An example is $g(z)=\sin (z)+\cos (z).\;$In Section 1
we completely characterize the RPE polynomials $p(z)$. Using a result from 
\cite{BG}, $p$\ must be either even, an odd polynomial plus a constant, or a
quadratic polynomial composed with an odd polynomial $s(z)$(Theorem \ref{T2}%
). For the latter case, a simple example is $p(z)=(z^{2}+2z)\circ
s(z),\;f(z)=\cos (2\pi \sqrt{z+1}).$\ Similar results are given for when $%
f\circ p$\ is odd(Theorem \ref{T3}).

Call an entire function $f(z)$\ left pseudo even(LPE) if $f\circ g$\ is even
for some entire function $g(z)$\ which is not even. Every even entire
function is LPE, but there are also LPE entire functions which are not even.
An example is the polynomial $f(z)=(z^{2}+2z)^{2}.\;$In Section 2 we discuss
the more general question of when $p\circ g$\ can be even, where $g(z)$\ is
entire and $p(z)$\ is a polynomial. However, we are not able to give a
complete characterization of the LPE polynomials, as we do for the RPE
polynomials. We can prove(Theorem \ref{BOREL}) that if $g\;$is a
transcendental entire function with finitely many zeroes, which is neither
even nor odd, then $p\circ g$ cannot be even for any non-constant polynomial 
$p.$ We can prove the simple result that $p(z)=z^{n},\;n$\ odd$,$\ cannot be
LPE(see Lemma \ref{B} ). It is interesting to note, however, that it is
possible for an odd entire function to be LPE. For example, $f(z)=\sin
\left( \dfrac{\pi }{2}z\right) $is LPE since $f\circ g$\ is even, where $%
g(z)=1-\sin \left( \dfrac{\pi }{2}z\right) .\;$

In Section 3 we extend some of our results to cyclic compositions. Let $N$\
be a prime number, and let $\omega $\ be a primitive $Nth$\ root of unity.
If $f(z)$\ is analytic at $0,$\ and $f(z)=z^{k}\sum_{j=0}^{\infty
}a_{j}z^{jN},$\ then $f$\ is called cyclic mod $N$. $\;f$\ is cyclic $\func{%
mod}2$\ if and only if $f$\ is even or odd. If $N\geq 3,\;f$\ is entire, and 
$p$\ is a polynomial, then necessary and sufficient conditions for $f\circ p$%
\ to be cyclic mod $N$\ are more restrictive than for the even or odd
case(Theorem \ref{3}).

In Section 4 we discuss even and odd compositions of rational functions,
where the results for polynomials do not extend in general. For example
there are rational functions $f(z)$\ and $g(z)$\ such that $%
f(0)=g(0)=0,\;f\circ g$\ is even, but neither $f$\ nor $g$\ are even. It is
indeed possible for $f\circ f$\ to be even when $f$\ is not even($f\left(
0\right) \neq 0$). We do not know if this is possible for the cyclic case in
general.

In Section 5 we discuss even compositions of polynomials in two variables.
Again, the results do not extend in general. There are polynomials $%
P(z,w),\;Q(z),$\ and $R(z)$, $P(0,0)=0,$\ $Q(0)=R(0)=0,$\ and $P,Q,$\ and $%
R\;$neither even nor odd, with $P(Q(z),R(z))$\ even. One can prove, however,
that if $Q(0,0)=0$\ and $P(Q(z,w))$\ is even, where $P$\ is a polynomial in
one variable, then $P$\ or $Q$\ must be even. Related questions along these
lines are homogeneous or symmetric compositions of polynomials in two or
more variables.

\section{Entire functions Composed with a Polynomial}

Let $f$\ be an entire function and $p$\ a polynomial. We shall prove that if 
$f\circ p$\ is even, then either $p$\ is even, odd, or a quadratic
polynomial composed with an odd polynomial. First we need the following
lemma.

\begin{lemma}
\label{L1}Suppose that $g$\ is analytic at $0,$ and $g$ is neither even, nor
odd plus a constant(i.e. $g^{\prime }$\ is neither even nor odd). Then there
are no constants $a$\ and $b$\ such that $g(z)=ag(-z)+b$\ for all $z$.
\end{lemma}

\proof%
. Write $g(z)=\sum_{k=0}^{\infty }\,a_{k}\,z^{k},$ and let $m$ and $n$ be
distinct positive integers with $m$ odd, $n$ even, and $a_{m}\neq 0\neq
a_{n}.$\ Then the expansion of $g(z)-ag(-z)$ about $0$ contains the terms $%
(1+a)z^{m}$ and $(1-a)z^{n},$ and hence no value of $a$ can make $%
g(z)-ag(-z) $ a constant. 
\endproof%
\smallskip

\begin{theorem}
\label{T2}Let $p$\ be a polynomial. Then there exists a non-constant entire
function $f$\ such that $f\circ p$\ is even if and only if one of the
following holds:
\end{theorem}

(A) $p$\ is even.

(B) $p(z)=s(z)+k,$\ where $s$\ is odd and $k$\ is a constant.

(C) $p(z)=\left( s(z)+d\right) ^{2}+k,$\ where $s$\ is odd and $d$\ and $k$\
are constants, $d\neq 0.$

Proof. We can assume $p\neq 0$, which is equivalent to (A), (B), and (C) all
holding. First we prove that if $f\circ p$\ is even, then one of (A), (B),
or (C) must hold. Assume first that neither (A) nor (B) holds, and let $%
q(z)=p(-z).$\ Then

\begin{equation}
f(p(z))=f(q(z))\ \forall z  \label{eqn1}
\end{equation}
By \cite{BG},

\begin{equation}
p(z)=a\,q(z)+b  \label{eqn2}
\end{equation}
or 
\begin{equation}
p(z)=r^{2}(z)+k,\,q(z)=\left( r(z)+c\right) ^{2}+k  \label{eqn3}
\end{equation}
By Lemma \ref{L1}, (\ref{eqn3}) must hold. Now $p(-z)=r^{2}(-z)+k=\left(
r(z)+c\right) ^{2}+k\Rightarrow r^{2}(-z)=\left( r(z)+c\right)
^{2}\Rightarrow r(-z)=\pm (r(z)+c)\Rightarrow -r^{\prime }(-z)=\pm r^{\prime
}(z)\Rightarrow r^{\prime }$\ is either even or odd. If $r^{\prime }$\ is
odd, then $r$\ is even, which implies that $p$\ is even. That contradicts
the assumption that (A) does not hold. Hence $r^{\prime }$\ is even $%
\Rightarrow r(z)=s(z)+d,$\ where $s$\ is an odd polynomial. Hence $%
p(z)=\left( s(z)+d\right) ^{2}+k,$\ and (C) holds. Now if (A) and (C) do not
hold, the above argument shows that $p$\ must satisfy (\ref{eqn2}). By Lemma 
\ref{L1}, $p(z)=s(z)+k,$\ where $s$\ is odd, and hence (B) holds. Finally,
if (B) and (C) do not hold, then Lemma \ref{L1} and the argument above show
that $p$\ must be even, and hence (A) holds.

To prove sufficiency, clearly if (A) or (B) hold, then one can always find a
non-constant entire function $f$\ such that $f\circ p$\ is even. If (C)
holds, one can choose $f(z)=\cos \left( \dfrac{2\pi }{d}\sqrt{z-k}\right) .$%
\ Then $f((z+d)^{2}+k)=\cos \left( \dfrac{2\pi z}{d}\right) $, which is
even. Hence $f((s(z)+d)^{2}+k)$\ is also even. 
\endproof%

\begin{remark}
(A) and (B) are equivalent to $p(z)-p(0)$\ being even or odd.
\end{remark}

\begin{remark}
Theorem \ref{T2} gives a complete characterization of the RPE polynomials.
\end{remark}

An odd version of Theorem \ref{T2} also follows. The proof is similar to the
proof of Theorem \ref{T2} and we omit it.

\begin{theorem}
\label{T3}Let $p$\ be a polynomial. Then there exists a non-constant entire
function $f$\ such that $f\circ p$\ is odd if and only if one of the
following holds:

(A) $p$\ is odd.

(B) $p(z)=\left( s(z)+d\right) ^{2}+k,$\ where $s$\ is odd and $d$\ and $k$\
are constants, $d\neq 0.$
\end{theorem}

\begin{remark}
If (B) of Theorem \ref{T3} holds, choose $f(z)=\cos \left( \dfrac{\pi }{2d}%
\sqrt{z-k}\right) $\ to make $f\circ p$\ odd.
\end{remark}

\begin{example}
Let $p(z)=z^{k}+z,$\ where $k\geq 4,\;k\;$even. Then it is not hard to show
that (C) of Theorem \ref{T2}(or (B) of Theorem \ref{T3}) does not hold(the
other conditions are trivial). Hence there is no non-constant entire
function $f$\ such that $f\circ p$\ is even or is odd. There is, however, a
function $f$\ analytic at $0$\ such that $f\circ p$\ is even(or odd). Since $%
p$\ has an inverse in a neighborhood of $z=0,$\ just let $f=x^{2}\circ
p^{-1} $(or $p^{-1}$).
\end{example}

\section{Polynomials Composed with an Entire Function}

Let $p$\ be a non-constant polynomial and $f$\ an entire function. As we
show below, $p\circ f$\ can be even, even when $f$\ is neither even, nor odd
plus a constant. It would be nice to have a theorem similar to Theorems \ref
{T2} and \ref{T3} which characterize when $p\circ f$\ is even, but we are
not able to do this. One can easily show that if $p\circ f$\ is even, then
there must be a polynomial relation between the even and odd parts of $f$,
but this is not sufficient. For example, let $f(z)=e^{z}=\sinh (z)+\cosh
(z). $\ There is a polynomial relation between $\sinh (z)$\ and $\cosh (z),$%
\ but no non-constant polynomial in $e^{z}$\ can be even. We prove a more
general result than this below. However, certain polynomial relations in the
even and odd parts of $f$\ do imply that $p\circ f$\ is even.

\begin{theorem}
\label{EO}Let $f(z)=E(z)+O(z),$\ where $E$\ is even, $O$\ is odd, and $%
E^{2}+O^{2}=1$\ for all $z.$\ Then $p\circ f$\ is even, where $%
p(z)=z^{4}-2z^{2}.$
\end{theorem}

\proof%
. $(E+O)^{2}-1=2OE,$\ which is odd. Hence$\left( (E+O)^{2}-1\right) ^{2}$is
even$\Rightarrow \left( (E+O)^{2}-1\right) ^{2}-1=(E+O)^{4}-2(E+O)^{2}\;$is
even as well. 
\endproof%

\begin{example}
Let $f(z)=\cos (z)+\sin (z).$\ Then $p(\sin (z)+\cos (z))=\allowbreak -4\cos
^{4}z+4\cos ^{2}z-1.\allowbreak $
\end{example}

The example above can be modified so that $f(0)=0$\ and $p$(as well as $f$)
is not even. Just let $f(z)=\cos (z)+\sin (z)-1,$\ and let$%
\,q(z)=z^{2}+2z.\, $Then$\,q(f(z))=2\cos z\sin z,$which is odd. Hence $%
p(f(z))$\ is even, where $p=q^{2}.$\ Note also that $f$\ is neither even,
nor odd plus a constant. However, $f$\ does have infinitely many zeros. This
must be the case, as the following theorem shows.

\begin{theorem}
\label{BOREL}Suppose that $f$ is a transcendental entire function with
finitely many zeroes in the plane, and let $p$ be a non-constant polynomial.
\end{theorem}

(A) If $f$ is neither even nor odd, then $p\circ f\;$

cannot be even.

(B) If $f$ is not odd, then $p\circ f\;$cannot be odd.

Before proving Theorem \ref{BOREL}, we need the following theorem first
stated by Borel. The first complete proof was given by R. Nevanlinna \cite
{NEV}.

\begin{theorem}
\label{BORE}(Borel \cite{BOR})Let $a_{i}(z)$ be an entire function of order $%
\rho ,$ let $g_{i}(z)$ also be entire and let $g_{j}(z)-g_{i}(z)(i\neq j)$
be a transcendental function or polynomial of degree higher than $\rho .$
Then
\end{theorem}

\[
\sum_{i=1}^{n}a_{i}(z)e^{g_{i}(z)}=a_{0}(z)
\]

holds only when

\[
a_{0}(z)=a_{1}(z)=\cdots =a_{n}(z)=0
\]

\proof%
of Theorem \ref{BOREL}. We prove (A), the proof of (B) following in a
similar fashion. Since $f$ is transcendental and has finitely many zeroes, $%
f(z)=Q(z)e^{g(z)},$ where $Q$ is a polynomial and $g$ is a non-constant
entire function. Write $p(z)=\sum_{i=0}^{n}a_{i}z^{i},\;a_{n}\neq 0,n\geq 1.$
Then

\[
p(f(z))=\sum_{i=0}^{n}a_{i}(Q(z))^{i}e^{ig(z)}
\]

Hence 
\begin{eqnarray*}
&&p(f(z))-p(f(-z)) \\
&=&\sum_{i=0}^{n}a_{i}(Q(z))^{i}e^{ig(z)}-%
\sum_{i=0}^{n}a_{i}(Q(-z))^{i}e^{ig(-z)}
\end{eqnarray*}

Now suppose that $p\circ f$ is even.

\underline{Case 1}: $g$ is even.

Then, since $p\circ f$ is even, 
\[
\sum_{i=0}^{n}a_{i}e^{ig(z)}((Q(z))^{i}-Q(-z)^{i})=0
\]
By Theorem \ref{BORE} with $a_{i}(z)=a_{i}(Q(z))^{i}-Q(-z)^{i},i\geq 1$ and $%
a_{0}=0$, which are entire functions of order $0$, $%
a_{i}(Q(z))^{i}-Q(-z)^{i}=0$ for all $i$. Since $a_{n}\neq
0,(Q(z))^{n}-Q(-z)^{n}=0,$ which implies that $Q$ is either even or odd(see
Lemma \ref{B} below). This contradicts the fact that $f$ is neither even nor
odd.

\underline{Case 2}: $g(z)=h(z)+C,$ where $h$ is odd and non-zero, and $C$ is
a constant.

Then, since $p\circ f$ is even, 
\[
\sum_{i=0}^{n}a_{i}e^{iC}(Q(z))^{i}e^{ih(z)}-%
\sum_{i=0}^{n}a_{i}e^{iC}(Q(-z))^{i}e^{-ih(z)}=0
\]
Again, by Theorem \ref{BORE}, with $a_{i}(z)=a_{i}e^{iC}(Q(z))^{i}$ or $%
a_{i}e^{iC}Q(-z)^{i},i\geq 1$ and $a_{0}=0$, $a_{i}e^{iC}(Q(z))^{i}=0$ for
all $i.$ Hence $(Q(z))^{n}=0,$ which implies that $Q=0,$ which contradicts
the fact that $f$ cannot be $0$.

\underline{Case 3}: $g$ is neither even, nor odd plus a constant.

By Lemma \ref{L1}, $kg(z)-jg(-z)$ cannot be a constant for $k$ and $j$
positive integers. Since $p\circ f$ even implies that

$\sum_{i=0}^{n}a_{i}(Q(z))^{i}e^{ig(z)}-%
\sum_{i=0}^{n}a_{i}(Q(-z))^{i}e^{ig(-z)}=0,$

by Theorem \ref{BORE}, with $a_{i}(z)=a_{i}(Q(z))^{i}$ or $%
a_{i}Q(-z)^{i},i\geq 1$ and $a_{0}=0$, $a_{i}(Q(z))^{i}=0$ for all $i.$
Hence $(Q(z))^{n}=0,$ which implies that $Q=0,$ which contradicts the fact
that $f$ cannot be $0$. 
\endproof%

\begin{remark}
Theorem \ref{BOREL} does \textit{not} follow in general if $f$ is not
transcendental, as the simple example $f(z)=z-1$\ and $p(z)=(z+1)^{2}$ or $%
p(z)=(z+1)^{3}$ shows. However, (A) does follow for polynomials with the
stronger assumption that $f$ is neither even nor odd plus a constant. (B)
follows for polynomials if $f$ is not odd plus a constant. This follows from 
\cite{HR} or from Theorem \ref{F}(B) in the next section.
\end{remark}

Finally we end this section by proving that if $p(z)=z^{n}$ and $f(z)$\ is
neither even nor odd, then $p\circ f$\ cannot be even.

\begin{lemma}
\label{B}Let $f$\ be an entire function, and suppose that $%
f^{\,n}(z)=(f(z))^{n}$\ is even for some positive integer $n$. Then $f$\
must be even or odd.
\end{lemma}

Proof.$(f(z))^{n}=(f(-z))^{n}\Rightarrow \left( \dfrac{f(z)}{f(-z)}\right)
^{n}=1\Rightarrow \dfrac{f(z)}{f(-z)}=e^{\tfrac{2\pi ki}{n}},$\ where $k\in
Z_{+}$\ depends on $z$, and $f(-z)\neq 0.$\ By continuity, $\dfrac{f(z)}{%
f(-z)}=e^{\tfrac{2\pi ki}{n}},$\ $k$\ independent of $z$. Hence $f(z)=e^{%
\tfrac{2\pi ki}{n}}f(-z).$\ Let $f(z)=E(z)+O(z),$\ where $E$\ is even and $O$%
\ is odd. Then $(1-e^{\tfrac{2\pi ki}{n}})\,E(z)=(-1-e^{\tfrac{2\pi ki}{n}%
})\,O(z).$\ If $e^{\tfrac{2\pi ki}{n}}\neq 1$\ or $-1,$\ then $E(z)=O(z)=0$\
and $f=0.$\ If $e^{\tfrac{2\pi ki}{n}}=1$, then $O(z)=0$\ and $f$\ is even.
If $e^{\tfrac{2\pi ki}{n}}=-1($\ $n$\ is even), then $E(z)=0$\ and $f$\ is
odd. 
\endproof%

\begin{remark}
Lemma \ref{B} says that $z^{n},\;n$\ odd, cannot be LPE.
\end{remark}

\section{\protect\medskip Cyclic Compositions}

Let $N$\ be a prime number, and let $\omega $\ be a primitive $Nth$\ root of
unity. Let

\[
C_{k}=\{\,f\;analytic\;at\;0:\ f(\omega z)=\omega
^{k}\,f(z)\,\},\,k=0,\ldots ,N-1
\]
$\,C=\bigcup {}_{k=0}^{N-1}\,C_{k}.\;$We call the functions in $C$\ cyclic
mod $N$. 
\[
f\in C_{k}\Rightarrow f(z)=z^{k}\sum_{j=0}^{\infty }a_{j}z^{jN}\,
\]
Of course, $C$\ and $C_{k}$\ depend on $N,$\ but we supress this in our
notation. For $N=2,\,C_{0}$\ are the even functions and$\,C_{1}$\ are the
odd functions, analytic at $0$. In \cite{HR} we proved that if $p$\ and $q$\
are polynomials, $q(0)=0,\,q$\ neither even nor odd, then $p\circ q\;$is
neither even(if $p\;$not constant) nor odd. As a corollary, if $q(0)=0$\ and 
$p\circ q\;$is even, then $p$\ and/or $q$\ must be even.

A similar result follows for the odd case. We extend this result now to
cyclic functions.

\begin{proposition}
\label{C}Let $p$\ and $q$\ be polynomials with $p$\ not constant, $q(0)=0,$\
and $q\notin C$. Then $p\circ q\notin C$.
\end{proposition}

Before proving Proposition \ref{C} we need the following lemma.

\begin{lemma}
\label{mod}Let $m,n,$\ and $r$\ be positive integers with $n\equiv 0\func{mod%
}N,\;r%
\not
%
\equiv m\func{mod}N$. Then $m(n-1)+r%
\not
%
\equiv 0\func{mod}N$.
\end{lemma}

Proof. $n=kN,\;r=lN+i,\;m=sN+j,$\ with $i\neq j$. Then$%
\,m(n-1)+r=(sN+j)(kN-1)+lN+i\equiv (i-j)\func{mod}N%
\not
%
\equiv 0\func{mod}p$. 
\endproof%

Proof of Proposition \ref{C}. Let $m=\deg q,$\ $n=\deg p$, and assume
w.l.o.g. that $p$\ is monic. Now write

\[
p(z)=(z-\alpha _{1})\cdots (z-\alpha _{n}),\;q(z)=\sum_{k=1}^{m}a_{k}\,z^{k}
\]

Then

\begin{equation}
p(q(z))={\LARGE (}\sum_{k=1}^{m}a_{k}\,z^{k}-\alpha _{1}{\LARGE )}\cdots 
{\LARGE (}\sum_{k=1}^{m}a_{k}\,z^{k}-\alpha _{n}{\LARGE )}  \label{E}
\end{equation}

We prove first that $p\circ q\notin C_{0}$.

\underline{Case 1} $m\equiv 0\func{mod}N$. Let $r$\ be the highest power of $%
z$\ in $q$, $r%
\not
%
\equiv 0\func{mod}N.\;$Then it follows easily from $(\ref{E})\;$that $%
m(n-1)+r$\ is the highest power of $z$\ in $p\circ q%
\not
%
\equiv 0\func{mod}N,$\ and the coefficient of $z^{m(n-1)+r}$\ is $%
n\,a_{m}^{n-1}a_{r}\neq 0$. Hence $p\circ q\notin C_{0}$.

\underline{Case 2} $m%
\not%
\equiv 0\func{mod}N$. Suppose that $p\circ q\in C_{0}$. We shall derive a
contradiction. $\deg (pq)=mn\equiv 0\func{mod}N\Rightarrow n\equiv 0\func{mod%
}N.\;$Let $r$\ be the highest power of $z$\ in $q$, $r%
\not
%
\equiv m\func{mod}N.$\ Note that $q(0)=0\Rightarrow r>0,$\ and thus there is
no cancellation of $z^{r}$\ with any $\alpha _{j}$\ in $(\ref{E}).$\ By
Lemma \ref{mod}, $m(n-1)+r%
\not
%
\equiv 0\func{mod}N$\ . Now $s>r\Rightarrow s\equiv m\func{mod}N.$\ Also, by 
$(\ref{E}),$\ any power of $z$\ in $p\circ q$\ larger than $m(n-1)+r$\ has
the form $s_{1}+\cdots +s_{n},$\ each $s_{j}>r.$\ Now suppose that $m\equiv k%
\func{mod}N.$\ Then $s_{1}+\cdots +s_{n}\equiv $

$n\,k\func{mod}N\equiv 0\func{mod}N.$\ Hence $m(n-1)+r$\ is the largest
power of $z$\ in $p\circ q$\ not congruent to $0\func{mod}N.$\ Since the
coefficient of $z^{m(n-1)+r}$\ is $n\,a_{m}^{n-1}a_{r}\neq 0$, $p\circ
q\notin C_{0}$. This is a contradiction, and hence $p\circ q\notin C_{0}$.

Now suppose that $p\circ q\in C_{j}$\ for some $j.$\ Then $p^{N}\circ
q=\left( p\circ q\right) ^{N}\in C_{0},$\ which contradicts the case just
proven. 
\endproof%

\begin{remark}
Proposition \ref{C} does not follow if $q(0)\neq 0$. For example, let $%
p(z)=(z+1)^{N},\,q(z)=z-1.$\ Also, it is trivial that Proposition \ref{C}
does not follow if $N$\ is not prime.
\end{remark}

\begin{remark}
The proof of Proposition \ref{C} shows that $q(0)=0$\ is only needed for
case 2, where $\deg (q)%
\not
%
\equiv 0\func{mod}N.$
\end{remark}

\begin{theorem}
\label{F}Let $p$\ and $q$\ be polynomials with $p\circ q\in C$.

(A) Suppose $q(0)=0.$\ Then $q\in C$\ or $p\in C.$\ Also, if $p$\ is not
constant and $q\notin C_{0},$\ then $p\in C$\ and $q\in C.$\ Finally, if $%
p\circ q\in C_{0}$, then $p\in C_{0}\;$or $q\in C_{0}.$

(B) Suppose $q(0)\neq 0.$\ If $\;p\;$is not constant and $q\notin C_{0}$,
then $q(z)=r(z)+q(0),p(z)=s(z-q(0))$\ where $r\in C,\,s\in C.$
\end{theorem}

Proof. To prove (A), if $p$\ is constant, then $p\in C.$\ If $p$\ is not
constant and $q\notin C,$\ then $p\circ q\notin C$\ by Proposition \ref{C}.
Hence $q\in C.$\ Now suppose that $p$\ is not constant. We just showed that $%
q\in C.$\ Now suppose that $q\in C_{j}$\ , $j\geq 1$\ and $p\circ q\in
C_{k}. $\ Then

\[
p(q(\omega z))=\omega ^{k}\,p(q(z))\ 
\]
and

\[
p(q(\omega z))=p(\omega ^{j}\,q(z))
\]
and hence $p(\omega ^{j}\,u)=\omega ^{k}\,p(u)$\ for any complex number $u.$%
\ Choose $r$\ so that $rj\equiv 1\func{mod}N.$\ Then 
\[
p(\omega u)=p(\omega ^{rj}\,u)=\omega ^{rk}\,p(u)
\]
which implies that $p\in C_{i},$\ where $rk\equiv i\func{mod}N.$\ Hence $%
p\in C.$\ Finally, suppose that $p\circ q\in C_{0}.$\ Again, $q\in C.$\ If $%
q\notin C_{0}$, then $q\in C_{j}$\ , $j\geq 1.$\ Then

\[
p(q(\omega z))=\,p(q(z))
\]
\ and 
\[
p(q(\omega z))=p(\omega ^{j}q(z))
\]
and hence $p(\omega u)=p(\omega ^{rj}u)=p(u),rj\equiv 1\func{mod}N.$\ Thus $%
p\in C_{0}.$

To prove (B), let $\;r(z)=q(z)-q(0)$. If $p$\ is not constant, then $%
p(q(z))=p(r(z)+q(0))=s(r(z)),$\ where$\;s(z)=p(z+q(0))$\ is a non-constant
polynomial. By (A), $r\in C$\ and $s\in C$\ since $r\notin C_{0}.$\ Hence $%
q(z)=$\ $r(z)+q(0)$\ and $p(z)=s(z-q(0))$\ with $r\in C$\ and $s\in C$. 
\endproof%

\begin{remark}
It is not true in general that if $p\circ q\in C_{j},$\ then $q\in C_{j}$\
or $p\in C_{j}.$\ This only holds in general for $j=0,$\ as (A) shows. For
the case $N=2,$\ however, it does follow that if $p\circ q$\ is odd and $%
q(0)=0$, then $p$\ or $q$\ must be odd(see \cite{HR}).
\end{remark}

\begin{theorem}
Let $p$\ be a polynomial with $p\circ p\in C.\;$

(A) If $p(0)=0$, then $p\in C.\;$

(B) If $p\circ p\in C_{0},$\ then $p\in C_{0}.$
\end{theorem}

Proof. (A) follows immediately from Theorem \ref{F}(A) with $p=q.$\ To prove
part (B), note that $\deg (p)\equiv 0\func{mod}N.$\ By Remark 2 following
the proof of Proposition \ref{C}(with $p=q$)$,\;p\in C_{0}.\;$%
\endproof%

\begin{remark}
(A) does not follow in general if $p(0)\neq 0.$\ For example, for $N=2$, $%
p(z)=-z+1$\ satisfies $p(p(z))=z$, which is odd, but $p$\ is neither odd nor
even.
\end{remark}

\subsection{Cyclic Compositions of Entire Functions with Polynomials}

For primes $N>2,$\ we now prove a theorem on cyclic compositions of entire
functions with polynomials $p(z)$. The possibilities for $p(z)$\ are more
restrictive than for Theorems \ref{T2} and \ref{T3}.

\begin{theorem}
\label{3}Let $p(z)\;$be a polynomial and $N\geq 3$\ a prime number. Then
there exists a non-constant entire function $f(z)$ such that $f\circ p\in C$%
\ if and only if $p(z)-p(0)\in C.$
\end{theorem}

Proof. $\Rightarrow )$\ Suppose first that $f\circ p\in C_{0}$ for some \
non-constant entire function $f(z),$ and write $p(z)=\sum_{k=0}^{n}\,a_{k}%
\,z^{k},\,a_{n}\neq 0.$\ Let $\omega $\ be a primitive $Nth$\ root of unity,
and let $q(z)=p(\omega z).$\ By \cite{BG} again, either $(\ref{eqn2})$\ or $(%
\ref{eqn3})$\ must hold. If $(\ref{eqn2})$\ holds, then

\[
p(\omega z)=\lambda \,p(z)+\beta
\]
\ which implies that

\[
\sum_{k=0}^{n}a_{k}\omega ^{k}z^{k}=\lambda
\sum_{k=0}^{n}\,a_{k}\,z^{k}+\beta
\]
\ Hence

\[
\lambda =\omega ^{n}\text{ and}\;a_{0}=\lambda a_{0}+\beta
\]
which implies that $\beta =a_{0}\left( 1-\omega ^{n}\right) $(Note that if $%
N\mid n,$\ then $\lambda =1,$\ $\beta =0,$\ and $p\in C_{0}$). Also,

$\lambda \,a_{j}=\omega ^{j}a_{j},\;j\geq 1\Rightarrow \omega ^{n}=\omega
^{j}$\ if $a_{j}\neq 0\Rightarrow n\equiv j\func{mod}N$\ if $a_{j}\neq
0,\;j\geq 1.$\ That implies that $p(z)-p(0)\in C_{k}$\ for some $k.$\ 

Now suppose that $(\ref{eqn3})$\ holds. Then 
\[
p(z)=r^{2}(z)+k,\,p(\omega z)=\left( r(z)+c\right) ^{2}+k
\]
\ for some constants $c$\ and $k.$\ Hence 
\[
r^{2}(\omega z)=(r(z)+c)^{2}
\]
which implies that 
\[
r(\omega z)=\pm (r(z)+c)
\]
and hence $\omega r^{\prime }(\omega z)=\pm r^{\prime }(z).$\ Letting $%
s(z)=r^{\prime }(z),$\ we have $\omega s(\omega z)=\pm s(z).$\ If $\omega
\,s(\omega z)=-s(z),$\ then $-\sum_{k=0}^{m}b_{k}z^{k}=\sum_{k=0}^{m}b_{k}%
\omega ^{k+1}z^{k},$\ where $s(z)=\sum_{k=0}^{m}b_{k}z^{k}.$\ This implies
that $-b_{j}=\omega ^{j+1}b_{j}.$\ Now $\omega ^{k}=-1$\ for some $k$
implies that 
\[
\cos \left( \dfrac{2\pi k}{N}\right) +i\sin \left( \dfrac{2\pi k}{N}\right)
=-1
\]
$\Rightarrow \dfrac{2\pi k}{N}$\ $=\pi l$\ for some integer $\smallskip
l\Rightarrow \cos \left( \dfrac{2\pi k}{N}\right) =\cos (\pi l)=1$\ since $%
N\geq 3$\ implies that $l$\ is even. Hence $\omega ^{k}\neq -1$\ for any $%
k\Rightarrow b_{j}=0$\ for all $j\Rightarrow r$\ is constant$\ \Rightarrow p$%
\ is constant $\Rightarrow p(z)-p(0)\in C_{0}.$\ 

Otherwise, 
\[
\omega s(\omega z)=-s(z)
\]
implies that $s(\omega z)=\omega ^{N-1}s(z)\Rightarrow s\in
C_{N-1}\Rightarrow r\in C_{0},$\ and hence $p($and$\;p(z)-p(0)\in C_{0}$\ as
well. If $f\circ p\in C_{k}$ for some $k\geq 1,$ then $f^{N}\circ p\in
C_{0}. $ By the case just proved, $p(z)-p(0)\in C_{k}$ for some $k$.

($\Leftarrow $\ Suppose that $p(z)-p(0)\in C_{k}$\ for some $k,$\ and let $%
\widetilde{p}(z)=p(z)-p(0).$\ \-Then

$\widetilde{p}(\omega z)=\omega ^{k}\,\widetilde{p}(z).$\ Let $%
f(z)=(z-p(0))^{N}.$\ Then $f(p(\omega z))=\left( \widetilde{p}(\omega
z)\right) ^{N}=\omega ^{kN}\left( \widetilde{p}(z)\right) ^{N}=\left( 
\widetilde{p}(z)\right) ^{N}$

$=(p(z)-p(0))^{N}=f(p(z))\Rightarrow $\ $f\circ p\in C_{0}.$

\begin{corollary}
Let $g(z)$ be an entire function which is periodic and cyclic $\func{mod}N$
for some prime number $N\geq 3.$ Then $g$ is constant.
\end{corollary}

Proof. Suppose that $g$ is not constant, and let $L$ be the period of $g.$
Let $f(z)=g(L\sqrt[N]{z})$ and $p(z)=(z+1)^{N}.$ Then $f\circ p=g(Lz)\in C,$
but $p(z)-p(0)\notin C,$ which contradicts Theorem \ref{3}. 
\endproof%

\begin{remark}
Much of Theorem \ref{F} follows from Theorem \ref{3}, but the latter theorem
requires deeper results from \cite{BG} than used in the proof of Theorem \ref
{F}.
\end{remark}

\section{Rational Functions}

By \cite[Theorem 1]{HR} (or by Theorem \ref{F} and the remark following, if $%
p$\ and $q$\ are polynomials with $q(0)=0$\ and $p\circ q$\ even(or odd),
then $p$\ or $q$\ must be even(or odd). This does not hold for rational
functions in general, however, as the following example shows.

\begin{example}
Let $f(z)=g(z)=\dfrac{z}{z-1}.$\ Then $f\circ g=z,$\ which is odd. Also, if $%
f(z)=\left( \dfrac{z}{z-1}\right) ^{2}$\ and $g(z)=\dfrac{z}{z-1},$\ then $%
f\circ g=z^{2},$\ which is even. In each case, $f(0)=g(0)=0,$\ and neither $f
$\ nor $g$\ are even or odd. For the cyclic case in general, just take $%
f(z)=\left( \dfrac{z}{z-1}\right) ^{p}$\ and $g(z)=\dfrac{z}{z-1}.$
\end{example}

The example above also shows that $f\circ f$\ can be odd even if $f(0)=0$\
and $f$\ is not odd. However, $f\circ f$\ cannot be even for any power
series(convergent or not) if $f(0)=0$(see \cite{REZ})$.$\ The following
example shows that $f\circ f$ can be even if $f(0)\neq 0.$

\begin{example}
Let $f(z)=\dfrac{z^{2}+z+1}{z^{2}-z+1}.$\ Then $f\circ f=\allowbreak \dfrac{%
3z^{4}+7z^{2}+3}{z^{4}+5z^{2}+1},$\ which is even, but $f$\ is not even.
Note that $f$\ and $f\circ f$\ are both analytic at $z=0$.
\end{example}

\section{Polynomials in Several Variables}

Let $P(z,w),$\ $Q(z),$\ and $R(z)$\ be polynomials with $P(0,0)=Q(0)=R(0)=0.$%
\ If $P(Q(z),R(z))$\ is even, must one of $P,\;Q,\;$or $R$\ be even ? The
answer is no, as is seen by the simple example $P(z,w)=z-w^{2},%
\;Q(z)=z^{6}+z^{4}+2z^{3}+z^{2},\;R(z)=z^{2}+z.$\ Then $P(Q(z),R(z))=z^{6},$%
\ which is even, but none of $P,Q,\;$or $R$\ are even. Note that the answer
is still no if we replace even throughout by even in each variable
separately. We can, however, prove the following result. We just prove the
two variable case, though the extension to any number of variables follows
easily. First we need the following lemma.

\begin{lemma}
\label{Lab}Suppose that $P(z,w)$\ is a polynomial such that $%
R_{a,b}(z)=P(az,bz)$\ is even for any complex constants $a$\ and $b.$\ Then $%
P(z,w)$\ is even.
\end{lemma}

Proof. Suppose that $P(z,w)$\ contains a homogeneous term of odd degree $%
k,\sum_{r+s=k}c_{r,s}z^{r}w^{s}.$\ Then the coefficient of $z^{r+s}$\ in $%
P(az,bz)$\ is

$\sum_{r+s=k}\,c_{r,s}\,a^{r}b^{s}.\;$Since $P(az,bz)$\ is even for any
complex constants $a$\ and $b,$\ the polynomial in $a$\ and $b,$\ $%
\sum_{r+s=k}\,c_{r,s}\,a^{r}b^{s},$\ is identically $0$. Then $c_{r,s}=0$\
for all $r$\ and $s$\ such that $r+s=k.$\ This implies that $P(z,w)$\ is
even. 
\endproof%

\begin{theorem}
\label{PQR}Let $R(z,w)=P(Q(z,w)),$\ where $P$\ is a polynomial in one
variable and $Q$\ is a polynomial in two variables. If $Q(0,0)=0$\ and $R$\
is even, then $P$\ or $Q$\ must be even.
\end{theorem}

Proof. Let $Q_{a,b}(z)=Q(az,bz)$\ and $R_{a,b}(z)=R(az,bz)=P(Q_{a,b}(z)).$\
Then $Q_{a,b}(0)=0$\ and $R_{a,b}(z)$\ is even, for any constants $a$\ and $%
b $. By \cite[Theorem 1]{HR} (or by Theorem \ref{F} and the remark
following), $P$\ or $Q_{a,b}$\ must be even. If $P$\ is even, we are
finished. If $P$\ is not even, then $Q_{a,b}$\ is even. By Lemma \ref{Lab}, $%
Q$\ must be even. 
\endproof%

\begin{remark}
One can also prove an odd version of Theorem \ref{PQR}, or more generally a
cyclic version, as we did for polynomials in one variable.
\end{remark}

\begin{remark}
Other related questions for polynomials in several variables are:

(A) If $R(z,w)=P(Q(z,w))$\ is homogeneous, $Q(0,0)=0,\;$prove that $P$\ or $Q
$\ must be homogeneous. This can be proved with techniques similar to those
used in this paper.

(B) If $R(z,w)=P(Q(z,w))$\ is symmetric, must $Q$\ be symmetric ?
\end{remark}

\section{Open Questions}

(1) Let $f(z)$\ be an entire function with $f\circ f$\ even, $f$\ not even.
If $f(0)\neq 0,$\ then $f$\ need not be even(see \cite{HR}). Must $%
f(z)=O(z)+c,$\ where $O$\ is odd and $c$\ is a constant ?

(2) Let $N\geq 3$\ be a given prime number. Does there exist a function $%
f(z)\notin C_{0},$analytic at $0,$such that $f\circ f\in C_{0}$ ?(If such an 
$f$\ exists, then by \cite{REZ} $f(0)\neq 0.)$ If yes, does there exist a
rational function $R(z)\notin C_{0},$\ such that $R\circ R\in C_{0}\;$?

(3) Discuss the questions in this paper when $(f\circ g)(z)=h(k(z)),\;$%
where\ $h$ and $k$\ are entire functions. When must $f$ or $g$ be an entire
function in $k$ ?

(4) Give a complete characterization of the RPE and LPE entire functions,
and in particular the LPE polynomials.

\bigskip


\begin{thebibliography}{9}
\bibitem{BG}  I. N. Baker and Fred Gross, `On factorizing entire functions',
Proc. London Math. Soc 17 (1968), 69-76.

\bibitem{BOR}  E. Borel, `Sur les z\'{e}ros des fonctions enti\`{e}res',
Acta Math. 20 (1897).

\bibitem{HR}  Alan L. Horwitz and Lee A. Rubel, `When is the composition of
two power series even?', J. Austral. Math. Soc. (Series A) 56(1994), 415-420.

\bibitem{NEV}  R. Nevanlinna, `Le th\'{e}or\`{e}me de Picard-Borel et la
th\'{e}orie des fonctions m\'{e}romorphes', Paris, 1929.

\bibitem{REZ}  Bruce Reznick, `When is the iterate of a formal power series
odd?', J,. Austral. Math. Soc. (Series A) 28 (1979), 62-66.
\end{thebibliography}
\end{document}